\documentstyle[11pt]{article}

\newcommand{\bfr}{\begin{flushright}}
\newcommand{\efr}{\end{flushright}}
\newcommand{\bfl}{\begin{flushleft}}
\newcommand{\efl}{\end{flushleft}}
\newcommand{\btm}{\begin{itemize}}
\newcommand{\etm}{\end{itemize}}
\newcommand{\be}{\begin{equation}}
\newcommand{\ee}{\end{equation}}
\newcommand{\bc}{\begin{center}}
\newcommand{\ec}{\end{center}}
\newcommand{\ba}{\begin{array}}
\newcommand{\ea}{\end{array}}
\newcommand{\bey}{\begin{eqnarray}}
\newcommand{\eey}{\end{eqnarray}}
\newcommand{\ben}{\begin{eqnarray*}}
\newcommand{\een}{\end{eqnarray*}}
\newcommand{\bt}{\begin{tabular}}
\newcommand{\et}{\end{tabular}}

\newcommand{\hm}{\hat m}

\newcommand{\hM}{\hat M}

\newcommand{\hDelta}{\hat \Delta}
\newcommand{\p}{{\cal P}}

\newcommand{\R}{R}
\newcommand{\A}{{\cal A}}
\newcommand{\EXP}{E}
\newcommand{\sign}{\mbox{sign}}

\begin{document}
\setcounter{page}{1}
\pagenumbering{roman}
 
\Large{\bf  Strongly-Consistent Nonparametric Forecasting and Regression
for Stationary Ergodic Sequences }

\vspace{0.3 in}
\Large{  Sidney Yakowitz, L\'{a}szl\'{o}  Gy\"orfi, John Kieffer and  Guszt\'{a}v Morvai }

\vspace{0.3 in}
\large{ J. Multivariate Anal.  71  (1999),  no. 1, 24--41.}

\vspace{1 in}

\bc {\bf Abstract} \ec

 Let $\{(X_i,Y_i)\}$ be a stationary ergodic time
series  with $(X,Y)$ values in the product space $\R^d\bigotimes \R .$
This study offers what is believed to be the first
 strongly consistent (with respect to
pointwise, least-squares, and uniform distance)
 algorithm for inferring $m(x)=E[Y_0|X_0=x]$ under the presumption 
that $m(x)$ is uniformly Lipschitz continuous.  
Auto-regression, or forecasting, is an important special case,
and as such our work extends the literature of nonparametric,
nonlinear forecasting by circumventing customary mixing assumptions.
The work is motivated by a time series model in stochastic finance and by
perspectives of its contribution to the issues of universal
time series estimation. 
\newpage
\pagenumbering{arabic}

\section{Introduction}

Nonparametric regression has been applied to a variety of 
contexts, in particular to time series modeling and prediction.
  The 
present study contributes to the methodology by 
showing how a regression function can be consistently inferred
from time series data under no  process assumptions
beyond stationarity and ergodicity.  (A Lipschitz condition on
the regression function itself will be imposed.)

Toward showing how our methodology can impinge on an 
 established  research area, we give one substantive
application  to a practical problem  in stochastic
finance:
Many works, such as  the Chapter 
 entitled ``Some Recent Developments in
Investment Research'' of the prominent text \cite{BKM}, argue for
the need to move beyond the Black-Scholes stochastic differential
equation. This and  other studies suggest
 the so-called ARCH and GARCH extensions
as a promising direction. The review of
this approach by Bollerslev  {\em et al.} \cite{BCK}  cites a litany of
unresolved issues. Of particular relevance is the discussion
of the need to account for persistency of the variance (Sections 2.6 and 3.6). 
(ARCH and GARCH models can be long-range dependent for certain 
ranges of parameters. In these cases, statistical analysis is
delicate \cite{davmik}.)

The basic idea behind the ARCH/GARCH setup is that one must
allow the asset volatility
(variance) to change dynamically, and perhaps (GARCH) to depend on current  
and past volatility values.  The review \cite{BCK} documents (p. 30) that
several authors have applied nonparametric and semiparametric
regression, with some success, to infer the 
ARCH functions from data. These methods
can fail if fairly stringent mixing conditions are not in force.
Masry and Tjostheim \cite{matj}, because of their rigorous
consideration of consistency, sets the stage for appreciating the potential of
the present investigation.
They propose that both the asset dynamics and volatility of a nonlinear
ARCH series be
inferred from nonparametric classes of
regression functions.  By imposing some fairly severe assumptions,
which would be tricky to validate from data, these authors are
able to assure that the ARCH process is strongly  mixing
(with exponentially decreasing parameter)
and consequently standard kernel techniques are applicable.

On another avenue toward asset series  modelling,
decades ago, Mandelbrot  suggested that
fractal processes  should be considered
in this context.
Fractals have been of interest to theorists and modellers alike
in part because they can display   persistency.
In his 1999 study, ``A Multifractal Walk down Wall Street,'' 
\cite{man}  Mandelbrot argues
 that conventional models for portfolio theory ignore soaring
volatility, and that is akin to a mariner ignoring the possibility
of a typhoon on the basis of the observation that weather is
moderate 95\% of the time.  

Such persistence as exhibited in the 
models of finance calls into question whether various processes of 
interest are
actually strongly mixing, a consistency requirement
 for conventional nonparametric
regression techniques.  We mention parenthetically that
telecommunications modelers are increasingly turning toward
long-range-dependent processes (e.g., \cite{res} and \cite{will})

As mentioned, the primary contribution  of the present paper is an algorithm
which is demonstratably consistent without imposition of
mixing assumptions.  The implication is that process assumptions
such as in \cite{matj} are not required for our algorithm.
The price paid for this flexibility is that convergence rates
and asymptotic normality cannot be assured.  This avenue is
worthy of exploration, nevertheless, because the limits of
process inference are clarified, and as a practical matter,
future work might lead to methods which are reasonably 
efficient if the process does satisfy  mixing assumptions,
but simultaneously assures convergence when mixing fails.

The algorithm is of the series-expansion type. 
The foundational idea (after Kieffer \cite{kief}) is that sometimes
it is possible to bound the error of ignoring the series tail,
and additionally  assure that the leading coefficients are
consistently estimated.   Specific constructs are given for
a partition-type estimator (Section 2) and for a kernel
series (Section 3).

We close this introduction with a survey of
the literature of nonparametric estimation for
stationary series without mixing hypotheses.

Let $ Y $ be a real-valued random variable  and let $ X $ 
be a $d$-dimensional random vector ({\em i.e.,}
the  observation or co-variate).
 We do not assume anything
about the distribution of $X$.  
As is customary in regression and forecasting,
 the main aim of the analysis here is to  minimize the mean-squared
error :
$$ \min_{f} E ( ( f(X) - Y)^{2} ) $$
over some space of real-valued functions $f(\cdot)$ defined on
the range of $X.$
This minimum is achieved by the
 regression function $ m(x) $, which
is  defined to be the conditional distribution of $Y$ given $X$:
\begin{equation}
 m(x) = E(Y \mid X=x ),\label{def}
\end{equation}
assuming  the expectation is well-defined, i.e., if $E |Y|<\infty.$
For each measurable function $ f $ one has
\begin{eqnarray*}
E( ( f(X) - Y)^{2}) & = & E ( ( m(X) - Y)^{2} )  + E (( m(X) - f(X) )^{2} ) \\
& = & E (( m(X) - Y )^{2} )  + \int (m(x) - f(x))^{2} \mu(dx),
\end{eqnarray*}
where $ \mu $ stands for the distribution of the observation $ X $. The
second term on the right hand side is called {\sl excess error} or 
{\sl integrated
squared error} for the function $f$, which is given the notation
\begin{equation}
J(f)=\int (m(x) - f(x))^{2} \mu(dx).\label{Jf}
\end{equation}
Clearly, the mean squared error for $f$ is close to that of the
optimal regression function  only
if the excess error
$J(f)$ is close to $0$.

With respect to the statistical problem of regression estimation,
 let $(X_1, Y_1 ),\ldots, (X_n, Y_n )\ldots$
be a stationary ergodic time series with marginal component denoted as  $(X,Y)$.
We study pointwise, $L_2(\mu ),$ and $L_{\infty}$ convergence
of the regression estimate $m_n$ to $m$. 
The estimator $m_n$ is called {\em weakly universally
consistent} if $J(m_n) \rightarrow  0$ in probability
for all distributions of $(X,Y)$ with $\EXP |Y|^2 < \infty$.
In the context of independent identically-distributed (i.i.d.)
 pairs $(X,Y),$
Stone \cite{st} first pointed out in 1977 that there exist
weakly universally consistent estimators.
Similarly, $m_n$ is called {\em strongly universally
consistent} if $J(m_n) \rightarrow  0$ a.s.
for all distributions of $(X,Y)$ with $\EXP |Y|^2 < \infty$.

  Following pioneering papers by Roussas \cite{rous69}
and Rosenblatt \cite{rosen},  a large
body of literature has accumulated on consistency and asymptotic
normality when the samples are correlated. 
In developments below, we will employ the notation, 
$$x_m^n=(x_m,\ldots,x_n),$$
presuming that $m\le n.$

The theory of nonparametric regression 
is of  significance in time series
analysis because, by 
considering  samples $\{( X_{n-q}^n,X_{n+1})\}$  in place of
the  pairs $\{(X_{n-q}^n,Y_n)\},$ 
the regression problem is transformed into
the  {\em forecasting} (or {\em auto-regression}) problem. 
 Thus, in forecasting,
 we are asking for the conditional
expectation of the next observation, given the $q-$past, with
$q$ a positive integer, or perhaps infinity.

As mentioned, nearly all the works on consistent statistical methods for forecasting
hypothesize {\sl mixing conditions,} which are assumptions about
how quickly dependency attenuates as a function of time separation
of the observables.
Under a variety of mixing assumptions, kernel and partitioning
 estimators 
are consistent, and have attractive rate properties.
The monograph by Gy\"orfi {\em et. al.} \cite{GyHSV89} gives a  
coverage of the literature of nonparametric inference for dependent
series. In that work, the partition estimator is shown to
be strongly consistent, provided $|Y|$ is a.s. bounded,
under $\phi-$ mixing and, with some provisos, under
$\alpha-$ mixing.
A drawback to much of the literature on nonparametric forecasting
 is that mixing conditions are  unverifiable
by available statistical procedures.   
Consequently, some investigators have examined the problem:

\vspace{0.4in}

\begin{quote}

  Let $\{X_i\}$ be a real vector-valued
 stationary ergodic sequence.  Find
a forecasting algorithm which is provably consistent in some sense.
\end{quote}
\vspace{0.4in}

Of course, some additional hypotheses regarding smoothness of
the auto-regression function and moment properties of the variables
will be allowed, but additional assumptions about attenuation of
 dependency are
ruled out. A {\em forecasting algorithm} for 
$$m(X^{-1}_{-p})=E[X_{0}|X^{-1}_{-p}]$$ here means
a rule giving a sequence $\{m_n\}$ of numbers such that 
for each $n$, $m_n$ is a measureable function 
determined entirely by the data segment
$X^{-1}_{-n}.$

For $X$ binary, Ornstein \cite{orn} provided a (complicated) 
strongly-consistent estimator
of\newline $E[X_0|X^{-1}_{-\infty}].$    Algoet \cite{bib1}  extended this
approach to achieve convergence over real-valued time series and in this 
and \cite{bib2}, connected the universal forecasting problem with 
fundamental issues in portfolio and gambling analysis as well as data
compression.
Morvai {\em et al.} \cite{MYGy96} offered another algorithm achieving
strong consistency in the above sense.  Their algorithm is
easy to describe and analyze, and such analysis shows, unfortunately,
that its data requirements make it infeasible \cite{mordis}.

On the negative side, Bailey \cite{Bailey} and Ryabko \cite{Ryabko88}
have proven that even over binary processes, there is no strongly
consistent estimator for the dynamic problem of inferring
$E[X_{n+1}|X^n_{0}]$, $n=0,1,2,\ldots.$

  We mention that for a real vector-valued 
Markov series with a stationary transition law,
 a strongly-consistent
estimator is available for inferring
$m(x)=E[X_0|X_{-1}=x]$  under the hypothesis that
the sequence is Harris recurrent \cite{yak1}.  
Admittedly this is a  dependency condition,
 but  the marginal (i.e., invariant) law 
need not exist:  Positive recurrence is not hypothesized. It is
difficult to imagine a  Markov condition weaker than Harris recurrence under
which statistical inference is assured.

It is to be noted that there are weakly-consistent estimators
for the moving regression problem $E[X_{n+1}|X^n_{0}]$,
 $n=0,1,2,\ldots$.
 It turns out that universal coding
algorithms (e.g. \cite{ziv}) of the information theory literature
can be converted to weakly-universally consistent algorithms
when the coordinate space is finite.  Morvai {\em et al} \cite{MYA}
have given a weakly-consistent (and potentially computationally
feasible) regression estimator  for 
the moving regression problem when  $X$ takes values from the set of real 
numbers. That work offers a synopsis of the literature of
weakly consistent estimation for stationary and ergodic
time series. All the studies we have cited on consistency without mixing
assumptions rely on algorithms which do not fall into any of
the traditional classes (partitioning, kernel, nearest neighbor) mentioned
in connection with i.i.d. regression.

From this point on,  
$\{(X_i,Y_i)\}$ will represent a time series with
  $(X,Y)$ values in $\R^d\bigotimes \R $
which is stationary and ergodic, and such that $E|Y_i|< \infty.$
 In Section 2,  we establish 
by means of a variation on
the partitioning method, that we have a.s. convergence pointwise,
and, in the case of bounded support, in  uniform distance,
provided that the regression function $m(x)=E[Y_0|X_0=x]$ satisfies a Lipschitz
condition and a bound on the Lipschitz constant is known in advance.  If 
furthermore $|Y|$ is known to be 
bounded (but perhaps the bound itself is not known), then our algorithm
converges in $L_2(\mu)$.   Section 3 provides analogous results
for a truncated kernel-type estimate. 
In summary, we miss our goal of
pointwise strong universal consistency only in that we must
restrict attention to regression functions satisfying a uniform
Lipschitz condition and the user must have a bound to the 
Lipschitz constant.
>From counter-examples in Gy\"orfi {\em et al.} \cite{gymy} one sees that some 
restrictions are needed. 

Recently 
we have obtained an important preprint by Nobel {\em et al.} \cite{NMK} which 
bears similarities with the present investigation. 
That study gives an algorithm for the long-standing problem of
density estimation of the marginal of a stationary sequence.
Somewhat analoguous to our conditions, Nobel {\em at al.} require that
the density function be of bounded variation.  The algorithm itself
is based on different principles from the present paper. 
In the paper \cite{MKN} by G.~Morvai, S.~Kulkarni, and A.~Nobel, the ideas 
in \cite{NMK} were  extended for regression estimation.

\section{Truncated partitioning estimation}

Let $(X_i,Y_i)_{i=1}^{\infty}$ be an  ergodic stationary
random sequence with $E|Y|<\infty$.   
Now we attack the problem of estimating the  regression function $m(x)$ 
by combining partitioning estimation with a series expansion.

Let $\p_k=\{A_{k,i} \ \mbox{i=1,\dots}\}$ be a nested cubic partition 
of $\R^d$  with volume $(2^{-k-2})^d$. 
Define  $A_k(x)$ to be the partition cell of $\p_k$  
into which $x$ falls.   Take
\bey
M_k(x) := E(Y|X\in A_k(x)).
\label{def_M}
\eey
One can show that
\bey
\label{Mkxtomx}
M_k(x)\to m(x)
\eey
for $\mu$-almost all $x\in \R^d$. (To see this, notice that
$\{M_k(X),\sigma(A_k(X)) \ k=1,2,\dots\}$
 is a martingale, $E|Y|<\infty$ implies $\sup_{k=1,2,\dots} E
|M_k(X)|<\infty$ and 
hence the martingale convergence theorem can be applied to achieve the
desired result (\ref{Mkxtomx}), cf. Ash \cite{ash72} pp. 292.)

For $k\ge 2$ let 
\bey
\label{def_delta}
\Delta_k(x)=M_k(x)-M_{k-1}(x).
\eey
Our analysis is motivated by the representation, 
\bey
\label{expansionofmx}
m(x)=M_1(x)+\sum_{k=2}^{\infty} \Delta_k(x) =\lim_{k\to\infty} M_k(x)
\label{rep}
\eey
 for $\mu$-almost all $x\in \R^d$. 
Now let $L>0$ be an arbitrary positive number. For integer $k\ge 2$ define
\bey
\Delta_{k,L}(x)=\sign(M_{k}(x)-M_{k-1}(x))\min(|M_{k}(x)-M_{k-1}(x)|,L2^{-k}). 
\eey
Define
\bey
m_{L}(x):=M_1(x)+\sum_{i=2}^{\infty} \Delta_{i,L}(x). 
\label{basic}
\eey
Notice that  $|\Delta_{i,L}(x)|\le L2^{-i},$ and hence $m_{L}(x)$ is 
well defined for all $x\in S$, where $S$ stands for the support of $\mu$ 
defined as  
\bey 
S:=\{x\in \R^d: \ \mu(A_k(x))>0 \ \mbox{for all $k\ge 1$.}\} 
\eey  
By Cover and Hart \cite{coha67}, $\mu(S)=1$.

The crux of the truncated partitioning estimate is inference of the terms 
$M_1(x)$ and $\Delta_{i,L}(x)$ for $i=2,3,\dots$ in
(\ref{basic}).  Define
\bey
\hM_{k,n}(x):={\sum_{j=1}^{n} Y_{j} 1_{\{X_j\in A_k(x)\}}
\over \sum_{j=1}^{n} 1_{\{X_j\in A_k(x)\}}}.
\eey
If $\sum_{j=1}^{n} 1_{\{X_j\in A_k(x)\}}=0,$ then take
$\hM_{k,n}(x)=0.$ 
Now for $k\ge 2$, define
\bey
\hDelta_{k,n,L}(x)=\sign(\hM_{k,n}(x)-\hM_{{k-1},n}(x))\min(|\hM_{k,n}(x)-\hM_{{k-1},n}(x)|,L2^{-k})
\eey
 and for $N_n$ a non-decreasing  unbounded sequence of positive integers,
define the estimator
\bey
\label{est_m}
\hm_{n,L}(x)=\hM_{1,n}(x)+\sum_{k=2}^{N_n} \hDelta_{k,n,L}(x).
\eey

\newtheorem{theorem}{Theorem}
\begin{theorem} \label{pointwisetheorem} Let $\{(X_i,Y_i)\}$ be a 
stationary ergodic time series 
with $E|Y_i|<\infty.$ Assume $N_n\to\infty$. 
Then almost surely,  
for all $x\in S$ 
\bey
\label{Lconvpoint}
\hm_{n,L}(x)\to m_{L}(x).
\eey
If the support $S$ of $\mu$ is a bounded subset of $\R^d$ then almost surely  
\bey
\label{Supconv}
\sup_{x\in S} |\hm_{n,L}(x)- m_{L}(x)|\to 0.
\eey
 If  either (i)  $|Y|\le D<\infty$ almost surely ($D$ 
need not be known)  or (ii) $\mu$ is of bounded support then 
\bey
\label{Lconvint}
\int (\hm_{n,L}(x)-m_{L}(x))^2 \mu(dx)\to 0. 
\eey
\end{theorem}
{\bf Proof }
First we prove that almost surely, for all $x\in S$, and for all $k\ge 1$,  
\bey
\label{hMkntoMk}
\lim_{n\to\infty} |\hM_{k,n}(x)- M_k(x)|=0.
\eey 
By the ergodic theorem, as $n\to \infty,$ a.s.,
$$
{\sum_{j=1}^{n} 1_{\{X_j\in A_{k,i}\}}\over n}\to P(X\in A_{k,i})=\mu(A_{k,i}).
$$  
Similarly,
$$
{\sum_{j=1}^{n} 1_{\{X_j\in A_{k,i}\}} Y_{j}\over n}\to
 E(Y 1_{\{X\in A_{k,i}\}})=\int_{A_{k,i}}m(z)\mu(dz),
$$ 
which is finite since $E|Y|$ is finite. 
Since there are countably many $A_{k,i}$,  
almost surely, for all 
$A_{k,i}\in \cup_{v} \p_v$ for which $\mu(A_{k,i})>0$:  
$$
{\sum_{j=1}^{n} 1_{\{X_j\in A_{k,i}\}} Y_{j}
 \over \sum_{j=1}^{n} 1_{\{X_j\in A_{k,i}\}} }
\to E(Y|X\in A_{k,i}).
$$
Since for each $x\in S$, $\mu(A_k(x))>0$ and for some index $i$, 
$A_k(x)=A_{k,i}$, we have proved (\ref{hMkntoMk}).

\noindent
Particularly, almost surely, for all $x\in S$, and for all $k\ge 2$,  
\bey
\label{hM1ntoM1}
\hM_{1,n}(x)\to M_1(x)
\eey
and 
\bey
\label{hDeltaknLtoDeltakL}
\hDelta_{k,n,L}(x)\to \Delta_{k,L}(x). 
\eey
Let integer $R>1$ be arbitrary. Let $n$ be so large that  $N_n>R$. For all $x\in S$, 
\begin{eqnarray}
\lefteqn{\nonumber |\hm_{n,L}(x)-m_{L}(x)|}\\
&\le&\nonumber  |\hM_{1,n}(x)-M_1(x)|
+\sum_{k=2}^{N_n} |\hDelta_{k,n,L}(x)-\Delta_{k,L}(x)|
+\sum_{k=N_n+1}^{\infty} |\Delta_{k,L}(x)| \\
&\le& \nonumber
 |\hM_{1,n}(x)-M_1(x)|
+\sum_{k=2}^{R} |\hDelta_{k,n,L}(x)-\Delta_{k,L}(x)|
+\sum_{k=R+1}^{\infty} (|\hDelta_{k,n,L}(x)|+|\Delta_{k,L}(x)|)\\
&\le&  \nonumber   |\hM_{1,n}(x)-M_1(x)|
+\sum_{k=2}^{R} |\hDelta_{k,n,L}(x)-\Delta_{k,L}(x)|
+ 2 L\sum_{k=R+1}^{\infty} 2^{-k}\\
&\le& \label{fundamentaleq}
 |\hM_{1,n}(x)-M_1(x)|
+\sum_{k=2}^{R} |\hDelta_{k,n,L}(x)-\Delta_{k,L}(x)| +L2^{-(R-1)}. 
\end{eqnarray}
By (\ref{hM1ntoM1}) and (\ref{hDeltaknLtoDeltakL}), almost surely, for all $x\in 
S$, 
\bey
 |\hM_{1,n}(x)-M_1(x)|
+\sum_{k=2}^{R} |\hDelta_{k,n,L}(x)-\Delta_{k,L}(x)|\to 0.
\eey
By (\ref{fundamentaleq}), almost surely, for all $x\in S$, 
\bey
 \limsup_{n\to\infty} |\hm_{n,L}(x)-m_{L}(x)|\le L2^{-(R-1)}. 
\eey
Since $R$ was arbitrary, (\ref{Lconvpoint}) is proved.  

Now we prove (\ref{Supconv}). Assume the support $S$ of $\mu$ is bounded. 
Let $\A_k$ denote the set of hyper-cubes  from partition $\p_k$           
with nonempty intersection with $S$. That is, define 
\bey
\A_k=\{A\in \p_k: \ A\cap S\neq\emptyset\}. 
\eey
Since $S$ is bounded, $\A_k$ is a finite set. For $A\in \p_k$ let 
$a(A)$ be the center of $A$. Then almost surely, 
\begin{eqnarray}
\lefteqn{ \nonumber \sup_{x\in S} \left( |\hM_{1,n}(x)-M_1(x)|+
\sum_{k=2}^{R} |\hDelta_{k,n,L}(x)-\Delta_{k,L}(x)|\right) }\\
&\le& \max_{A\in \A_1} |\hM_{1,n}(a(A))-M_1(a(A))|
+\sum_{k=2}^{R} \max_{A\in \A_k} |\hDelta_{k,n,L}(a(A))-\Delta_{k,L}(a(A))|\\
&\to& 0 
\end{eqnarray}
keeping in mind that only finitely many terms are involved in the 
maximization operation. The rest of the proof  goes virtually as before. 

Now we prove (\ref{Lconvint}). 
$$
|\hm_{n,L}(x)-m_L(x)|^2\le 2\left(|\hM_{1,n}(x)-M_1(x)|^2+
|M_1(x)+\sum_{k=2}^{N_n}\hDelta_{k,n,L}(x)-m_L(x)|^2\right).
$$
If condition (i) holds, then   for the first term we have dominated convergence
$$
 |\hM_{1,n}(x)-M_1(x)|^2\le (2D)^2,
$$
and for the second one, too: 
\begin{eqnarray*}
\lefteqn{ 
|M_1(x)+\sum_{k=2}^{N_n}\hDelta_{k,n,L}(x)-m_L(x)|}\\
&\le&\sum_{k=2}^{\infty} (|\hDelta_{k,n,L}(x)|+|\Delta_{k,L}(x)|)\\
&\le& L,
\end{eqnarray*}
and thus (\ref{Lconvint}) follows by Lebesgue's dominated convergence theorem,
$$
0=\int \lim_{n\to\infty} |\hm_{n,L}(x)-m_L(x)|^2 \mu(dx)= \lim_{n\to\infty}
\int |\hm_{n,L}(x)-m_L(x)|^2 \mu(dx)
$$
almost surely.   
If condition (ii) holds then (\ref{Lconvint}) follows from   (\ref{Supconv}). 

\noindent$\Box$

\newtheorem{corollary}{Corollary}
\begin{corollary} \label{pointwisecorollary} 
Assume $m(x)$ is Lipschitz continuous with Lipschitz constant 
$C$. 
With the choice of $L\ge C\sqrt{d}$,  for all $x\in S$, $m_L(x)=m(x)$ and 
Theorem ~\ref{pointwisetheorem} 
holds with $m_L(x)$ replaced by $m(x)$.  \end{corollary}
{\bf Proof } 
 Since $m(x)$ is Lipschitz with constant $L/\sqrt{d}$, for $x\in S$, 
\begin{eqnarray*}
|M_k(x)-m(x)|&\le& |{\int_{A_k(x)} m(y) \mu(dy)\over \mu(A_k(x))} -m(x)|\\
&\le& {1\over \mu(A_k(x))} \int_{A_k(x)} |m(y)-m(x)| \mu(dy)\\
&\le& {1\over \mu(A_k(x))} \int_{A_k(x)} (L/\sqrt{d})(2^{-k-2}\sqrt{d}) \mu(dy)\\
&=& L2^{-k-2}
\end{eqnarray*}
and $M_k(x)\to m(x)$. For $x\in S$ we get 
\begin{eqnarray*}
|M_k(x)-M_{k-1}(x)|&\le& 
|M_k(x)-m(x)|+ |m(x)-M_{k-1}(x)| \\
 &\le& L2^{-k-2}+L2^{-k-1}\\
&<& L2^{-k}.
\end{eqnarray*}
Thus $m(x)=M_1(x)+\sum_{k=2}^{\infty} \Delta_k(x) $ and 
$\Delta_{k,L}(x)=\Delta_k(x)$ for all $x\in S $. 
Hence for all $x\in S$, 
$$m_L(x)=
M_1(x)+\sum_{k=2}^{\infty} \Delta_{k,L}(x)= 
M_1(x)+\sum_{k=2}^{\infty} \Delta_k(x)=m(x)
$$ 
and Corollary \ref{pointwisecorollary} is proved. 

\noindent$\Box$

\noindent
{\bf Remark 1.}
If there is no truncation, that is if $L=\infty$, then $\hm_n=\hM_{N_n,n}$.
In this case, 
$\hm_n$ is the standard partitioning estimate (defined, for example
in \cite{GyHSV89}).
It is known that there is  an  ergodic process $(X_i,Y_i)$ with Lipschitz  
continuous $m(x)$ with constant $C=1$ such that a classical partitioning estimate 
is not even weakly  consistent. (cf. Gy\"orfi, Morvai, Yakowitz 
\cite{gymy}).

\noindent
{\bf Remark 2.}
Our consistency is not universal, however, since $m$ 
is hypothesized to be Lipschitz
continuous.

\noindent
{\bf Remark 3.}
$N_n$ can be data dependent, provided $N_n\to \infty$ a.s.

\noindent
{\bf Remark 4.}
The methodology here is applicable to linear auto-regressive
processes. 
Let $\{Z_i\}$ be i.i.d. random variables with $EZ=0$ and 
$Var(Z)<\infty$. Define 
\be
\label{autoregressiveeq}
W_{n+1}=a_1 W_n+a_2 W_{n-1}+\dots +a_K W_{n-K+1}+Z_{n+1}
\ee
where $\sum_{i=1}^{K} |a_i|< 1.$
Equation (\ref{autoregressiveeq}) yields a stationary ergodic solution. 
Assume 
$K\le d$. Let $Y_{n+1}=W_{n+1}$, and $X_{n+1}=(W_n,\dots,W_{n-d+1})$. 
Now 
\begin{eqnarray*}
m(X_{n+1})=E(Y_{n+1}|X_{n+1})&=&E(W_{n+1}|W_n,\dots,W_{n-d+1})\\
&=& a_1 W_n+a_2 W_{n-1}+\dots a_K W_{n-K+1}.
\end{eqnarray*}
The regression function $m(x)$ is Lipschitz continuous with constant 
$C=1$, since for $x=(x_1,\dots,x_d)$ and $z=(z_1,\dots,z_d)$, 
$$
|m(x)-m(z)|\le \sum_{i=1}^{K} |a_i|  |x_i-z_i|\le \max_{1\le i\le 
d}|x_i-z_i|\le \|x-z\|. 
$$
 
\section{Truncated kernel estimation}

Let $K(x)$ be a non-negative continuous kernel function 
with \[
b 1_{\{x\in S_{0,r}\}}\le K(x)\le  1_{\{x\in S_{0,1}\}},
\]  
where $0<b\le 1$ and $0<r< 1$. ( $S_{z,r}$ denotes the closed 
ball around $z$ with radius $r$.) 

\noindent
Choose \[
h_k=2^{-k-2}
\]
and
\bey
M_k^*(x)={E(YK({X-x\over h_k}))\over E(K({X-x\over h_k}))}
={\int m(z)K({z-x\over h_k})\mu(dz)\over \int K({z-x\over h_k})\mu(dz)}.
\eey
Let
\bey
\label{def_delta_star}
\Delta_k^*(x)=M_k^*(x)-M_{k-1}^*(x).
\eey
As a motivation, we note that Devroye \cite{de81} yields  
(\ref{Mkxtomx}), and therefore (\ref{expansionofmx}), too. \noindent
Now for $k\ge 2$, define
\bey
\Delta_{k,L}^*(x)=\sign(M_{k}^*(x)-M_{k-1}^*(x))\min(|M_{k}^*(x)-M_{k-1}^*(x)|,L2^{-k}). 
\eey
Define 
\bey
m_{L}^*(x):=M_1^*(x)+\sum_{i=2}^{\infty} \Delta_{i,L}^*(x). 
\eey
Put
$$\hM_{k,n}^*(x):={\sum_{j=1}^{n} Y_{j} K({X_j-x\over h_k})
\over \sum_{j=1}^{n} K({X_j-x\over h_k})}$$
where we use the convention that $0/0=0$. 
Now for $k\ge 2$, introduce 
\bey
\hDelta_{k,n,L}^*(x)=\sign(\hM_{k,n}^*(x)-\hM_{{k-1},n}^*(x))\min(|\hM_{k,n}^*(x)-\hM_{{k-1},n}^*(x)|,L2^{-k})
\eey
 and 
\bey
\label{est_m_star}
\hm_{n,L}^*(x)=\hM_{1,n}^*(x)+\sum_{k=2}^{N_n} \hDelta_{k,n,L}^*(x).
\eey
Redefine the support $S$ of $\mu$ as 
\bey 
S:=\{x\in \R^d: \ \mu(S_{x,1/k})>0 \ \mbox{for all $k\ge 1$}\}. 
\eey  
By Cover and Hart \cite{coha67}, $\mu(S)=1$.

\begin{theorem} \label{kernel} Let $\{(X_i,Y_i)\}$ be a 
stationary ergodic time series with $E|Y_i|<\infty.$ 
Assume $N_n\to\infty$.  Then   almost surely, 
for all $x\in S$,
\bey
\label{Lconvpoint2}
\hm_{n,L}^*(x)\to m^*_{L}(x).
\eey
If the support $S$ of $\mu$ is a bounded subset of $\R^d$ then almost surely
\bey
\label{sup*}
\sup_{x\in S} |\hm_{n,L}^*(x)-m_{L}^*(x)|\to 0.
\eey
 If either (i) $|Y|\le D<\infty$ almost surely ($D$ need not be 
known)  or (ii) $\mu$ is of bounded support then 
\bey
\label{Lconvint2}
\int (\hm_{n,L}^*(x)-m^*_{L}(x))^2 \mu(dx)\to 0. 
\eey
\end{theorem}
{\bf Proof } We first prove that~(\ref{hMkntoMk}) holds with $\hM_{k,n}^*$ 
and $M_{k}^*$.
Let 
$$
g_{k,n}(x)= {1\over n} \sum_{j=1}^{n}  Y_j K\left({X_j-x\over h_k}\right)
$$
and 
$$
g_k(x) =E \left( YK\left({X-x\over h_k}\right)\right).
$$
Similarly put 
$$
f_{k,n}(x)={1\over n} \sum_{j=1}^n K\left({X_j-x\over h_k}\right)
$$
and
$$
f_k(x)=E K\left({X-x\over h_k}\right).
$$
We have to show that almost surely, for all $k\ge 1$, and for all 
$x\in S$, both $g_{k,n}(x)\to g_k(x)$ and $f_{k,n}(x)\to f_k(x)$. 
Consider $g_{k,n}(x)$ with $k$ fixed. Let $Q\subseteq \R^d$ denote the set of 
vectors with rational coordinates. (Note that the set $Q$ has countably
many elements.) By the ergodic theorem, almost surely, for all  $r\in Q$, 
$$
g_{k,n}(r)\to g_k(r).
$$
Let $\delta>0$ be arbitrary. 
Let integers $Z-1>M>0$ be so large that  
$E\left( |Y| 1_{\{X\notin S_{0,M}\}}\right) <\delta$. 
By ergodicity, almost surely, 
$$
\sup_{x\notin S_{0,Z}} |g_{k,n}(x)|\le 
{1\over n} \sum_{i=1}^{n} |Y_i| 1_{\{X_i\notin S_{0,M}\}}\\
\to   E\left( |Y| 1_{\{X\notin S_{0,M}\}}\right)  < \delta.
$$
Since $K_{h_k}(x)=K({x\over h_k})$ is 
continuous and $K_{h_k}(x)=0$ if $\|x\|>h_k$ and  hence $K_{h_k}(x)$ is 
uniformly  continuous on $\R^d$. 
Define 
$$
U_k(u)=\sup_{x,z\in \R^d:\|x-z\|\le u} |K_{h_k}(x)-K_{h_k}(z)|.
$$ 
Let $B_{\delta}\subseteq S_{0,Z}\cap Q$ be a finite subset of vectors 
with rational coordinates  
such that 
$$
\sup_{x\in S_{0,Z}} 
\min_{r\in B_{\delta}} U_k(\left\|x-r\right\|)<\delta.
$$
For $x\in S_{0,Z}$, let $r(x)$ denote one of the closest rational vector 
$r\in B_{\delta}$ to $x$.   Now 
\begin{eqnarray*}
\sup_{x\in S_{0,Z}}|g_{k,n}(x)-g_{k,m}(x)|&\le& 
\sup_{x\in S_{0,Z}}  | g_{k,n}(x) -g_{k,n}(r(x))| \\
&+&
\sup_{x\in S_{0,Z}}  |g_{k,n}(r(x))-g_{k,m}(r(x))|\\
&+&\sup_{x\in S_{0,Z}}  |g_{k,m}(r(x))-g_{k,m}(x)|\\
&\le& \delta {1\over n} \sum_{i=1}^n |Y_i|  + 
\max_{r\in B_{\delta}} |g_{k,m}(r)-g_{k,n}(r)|+
\delta {1\over m} \sum_{i=1}^n |Y_i| .
\end{eqnarray*}
Combining the results, by the ergodic theorem, for almost all 
$\omega\in\Omega$, 
there exists $N(\omega)$ such that for all $m>N$, and $n>N$,  
\begin{eqnarray*}
\sup_{x\in \R^d}   |g_{k,n}(x)-g_{k,m}(x)|&\le& 
\sup_{x\in S_{0,Z}}   |g_{k,n}(x)-g_{k,m}(x)|\\
&+& 
\sup_{x\notin S_{0,Z}}   |g_{k,n}(x)-g_{k,m}(x)|\\
&\le&
2\delta E|Y|+ 3\delta.
\end{eqnarray*}
Since $\delta$ was arbitrary, for almost all $\omega\in\Omega$, for every 
$\epsilon>0$, there exists an integer $N_{\epsilon}(\omega)$ such that for 
all $m>N_{\epsilon}(\omega)$, $n>N_{\epsilon}(\omega)$:  
\bey
\label{uniformcauchygknm}
\sup_{x\in \R^d}   |g_{k,n}(x)-g_{k,m}(x)|<\epsilon.
\eey
As a consequence, almost surely, the sequence of functions 
$\{g_{k,n}\}_{n=1}^{\infty}$ 
converges uniformly. Since all $g_{k,n}$ are continuous, 
the limit function must be also continuous. Since almost surely, for 
all  $r\in Q$, 
$g_{k,n}(r)\to g_k(r)$, and by the Lebesgue dominated convergence $g_k$ 
is continuous,  the limit function must be $g_k$. 
Since there are countably many $k$, almost surely, for all $k\ge 1$, 
$$
 \sup_{x\in \R^d}| g_{k,n}(x)- g_k(x)|\to 0.
$$
The same  argument implies that almost surely, for all $k\ge 1$, 
$$
\sup_{x\in \R^d} |f_{k,n}(x)- f_k(x)|\to 0. 
$$  
We have proved (\ref{hMkntoMk}).  
The rest of the proof of ~(\ref{Lconvpoint2}) goes as in the proof of 
Theorem 1.
Now we prove ~(\ref{sup*}). 
Since now, by assumption, the support is bounded, and since it is  
closed, and hence it is compact. 
Now note that 
there must exist an $\epsilon>0$ such that $\inf_{x\in S} 
f_k(x)>\epsilon$. (Otherwise, there would be 
 a sequence $x_i\in S$ such that 
$\liminf_{i\to\infty} f_k(x_i)=0$. Continuity on a compact set would imply 
that there would be an $x\in S$ such that $f_k(x)=0$  in contradiction to 
the hypothesis that $x\in S$. ) By uniform convergence,  for large $n$,  $\inf_{x\in S} 
f_{k,n}(x)>\epsilon/2$. Thus  
\begin{eqnarray*}
\sup_{x\in S}
 \left| { g_{k,n}(x)\over f_{k,n}(x) }-{ g_k(x)\over f_k(x) }\right|
&\le&\sup_{x\in S} \left| { g_{k,n}(x) (f_{k}(x)/ f_{k,n}(x) )-g_k(x)\over 
f_k(x) }\right|\\
&\le& {1\over \epsilon}
\sup_{x\in S} \left| { f_{k}(x)\over f_{k,n}(x) }\right| 
|g_{k,n}(x)-g_k(x)|+|g_k(x)| \left| {f_{k}(x)\over f_{k,n}(x)}-1\right|\\
&\le&
{2 \over \epsilon^2} \sup_{x\in S}  |g_{k,n}(x)-g_k(x)|+
\sup_{x\in S} |g_k(x)| {2\over \epsilon} \sup_{x\in S} |f_{k,n}(x)-f_k(x)|\\
 &\to& 0.
 \end{eqnarray*}
Thus almost surely, for all $k\ge 1$, 
$$
\sup_{x\in S} |\hM_{k,n}^*(x)-\hM_{k}^*(x)|\to 0.
$$
Almost surely, for arbitrary integer $R>2$, 
$$
\sup_{x\in S} \left( |\hM_{1,n}^*(x)-M_1^*(x)|+
\sum_{k=2}^{R} |\hDelta_{k,n,L}^*(x)-\Delta_{k,L}^*(x)|\right)\to 0.
$$ 
The rest of the proof goes exactly as in Theorem 1.   

\noindent$\Box$

\begin{corollary} \label{kernelcor} 
Assume $m(x)$ is Lipschitz continuous with Lipschitz constant $C$.  
With the choice of $L\ge C$ for  all $x\in S$, $m_L^*(x)=m(x)$ and 
Theorem \ref{kernel} holds with $m_L^*(x)$ substituted by $m(x)$.
  \end{corollary}
{\bf Proof } 
 Since $m(x)$ is Lipschitz with constant $C$, for $x\in S$,  
\begin{eqnarray*}
|M_k^*(x)-m(x)|&\le& |{\int m(z)K({z-x\over h_k})\mu(dz)
\over \int K({z-x\over h_k})\mu(dz)} -m(x)|\\
&\le& {\int |m(z)-m(x)|K({z-x\over h_k})\mu(dz)
\over \int K({z-x\over h_k})\mu(dz)}\\
&\le&  C h_k\\
&\le& L 2^{-k-2},
\end{eqnarray*}
therefore
\[
|M_k^*(x)-M_{k-1}^*(x)|< L2^{-k}.
\]
The rest of the proof goes as in Corollary \ref{pointwisecorollary}.
 
\noindent$\Box$

\section{Conclusions}

This contribution is part of a long-standing endeavor of
the authors to extend nonparametric forecasting methodology
to the most lenient assumptions possible.
The present work does push into new territory: 
strong consistency for finite regression 
under a Lipschitz assumption.  The computational aspects
have not been explored, but the algorithms are so close
to their traditional partitioning and kernel counterparts
that it is evident that they could be implemented and
in fact, might be competitive.

The fundamental formula (\ref{basic}) leading to the
truncated histogram approach was motivated by a representation used
in a related but non-constructive setting by Kieffer \cite{kief}.
The essence is to see that an infinite-dimensional nonparametric
space may sometimes be decomposed into sums of terms in  finite dimensional
spaces, with tails of the summations being {\em a priori} 
asymptotically bounded over the regression class of interest.
Through different devices, two ideas for obtaining such tail bounds
for the partition and kernel methods have been presented.

Our contribution has been to apply the idea with Lipschitz
continuity assuring the negligibility. Thus, results here
are fundamentally intertwined with the Lipschitz bounds.  Perhaps
other useful expansions are possible.  The interplay of
finite subspaces and {\em a priori} bounded  tails has proven
a bit delicate.  Sections 2 and 3 present different attacks
to the error-bounding problem. 
The obvious nearest-neighbor estimator did not
yield to this technique because the radii are random and do not necessarily
decrease rapidly enough to assure bounded tails.  
The device which was successful here may find other applications.
Evidently, a similar investigation could be carried out for
regression classes having
Fourier expansions with coefficients vanishing sufficiently quickly.

It is well-known (e.g., \cite{Sh}) that universal convergence rates under the
generality of mere ergodicity do not exist.  An avenue which
would be worth exploring is that of adapting universal algorithms,
such as explored and referenced here, so that they
asymptotically attain the fastest possible convergence if, unknown to the
statistician, the time series happens to fall into a mixing class.
The design should be such that
 consistency is still assured if mixing rates do not hold.

\end{document}